\documentclass[reqno,10pt]{amsart}
\usepackage[utf8]{inputenc}
\usepackage{amsmath, amssymb, amsthm,amsfonts} 
\usepackage[english]{babel}
\usepackage{bbm,bm}
\usepackage{graphicx}
\usepackage{url}
\usepackage{epstopdf}
\usepackage[a4paper,bindingoffset=0.2cm,left=1.2cm,right=1.2cm,top=2.5cm,bottom=2cm,footskip=.8cm]{geometry}
\usepackage{subcaption}
\usepackage{textcomp}
\usepackage[]{algorithm2e}
\usepackage{multicol}
\usepackage{multirow}
\usepackage{wrapfig}
\usepackage{diagbox}
\usepackage{floatrow}
 \usepackage[comma]{natbib}
 
 \usepackage[titletoc,toc,title]{appendix}
\newfloatcommand{capbtabbox}{table}[][\FBwidth]

\usepackage{tikz}
\usetikzlibrary{arrows, automata,positioning,calc,shapes,decorations.pathreplacing,decorations.markings,shapes.misc,petri,topaths}
\usepackage{pgfplots}
  \pgfplotsset{compat=newest}
  \usetikzlibrary{plotmarks}
  \usepackage{grffile}
\newlength\figureheight
  \newlength\figurewidth
\pgfplotsset{%
    tick label style={font=\scriptsize},
    label style={font=\footnotesize},
    legend style={font=\footnotesize},
         every axis plot/.append style={very thick}
}
\usepackage{rotating}
\usepackage{amsbsy,enumerate}
\usepackage{graphicx}
\usepackage{comment}
\usepackage{mathrsfs} 
\usepackage{tikzscale}
\usepackage{xcolor}

\newcommand{\vb}{\vspace{5mm}}

\newcommand\myeq{\mathrel{\overset{\makebox[0pt]{\mbox{\normalfont\tiny\sffamily d}}}{=}}}

\allowdisplaybreaks

\makeatletter
\newtheorem*{rep@theorem}{\rep@title}
\newcommand{\newreptheorem}[2]{%
\newenvironment{rep#1}[1]{%
 \def\rep@title{#2 \ref{##1}}%
 \begin{rep@theorem}}%
 {\end{rep@theorem}}}
\makeatother

\makeatletter
\renewcommand*\env@matrix[1][\arraystretch]{%
  \edef\arraystretch{#1}%
  \hskip -\arraycolsep
  \let\@ifnextchar\new@ifnextchar
  \array{*\c@MaxMatrixCols c}}
\makeatother

\newtheorem{theorem}{Theorem}
\newreptheorem{theorem}{Theorem}

\newtheorem{remark}{Remark}
\newtheorem*{example}{Example}

\newreptheorem{proposition}{Proposition}

\begin{document}

\title[Input Rate Control in Stochastic Road Traffic Networks: Effective
Bandwidths]
{Input Rate Control in Stochastic Road Traffic Networks: \\
Effective Bandwidths}

\author{Nikki Levering, Sindo N\'u\~nez Queija}

\maketitle

\begin{abstract}
In road traffic networks, large traffic volumes may lead to extreme delays. These severe delays are caused by
the fact that, whenever the maximum capacity of a road
is approached, speeds drop rapidly. 
Therefore,  
the focus in this paper is on real-time \textit{control} 
of traffic input rates, thereby aiming to prevent such detrimental capacity drops.
To account for the fact that, by the heterogeneity
within and between traffic streams, 
the available capacity of a road suffers from randomness,
we introduce a stochastic flow model that describes
the impact of traffic input streams on the available road capacities.
Then, exploiting similarities with traffic control
of telecommunication networks, in which the available
bandwidth is a stochastic function of the input rate,
and in which the use of \textit{effective bandwidths} have proven an
effective input rate control framework,
we propose a similar traffic rate control policy
based on the concept of effective bandwidths.
This policy allows for increased waiting times at the 
access boundaries of the network, 
so as to limit the 
probability of large delays within the network.
Numerical examples show that, by applying such a control
policy capacity violations are indeed rare, and 
that the increased waiting at the boundaries of the network 
is of much smaller scale, compared to uncontrolled delays in the network.

\vb


\vb

\noindent
{\sc Affiliations.} 
{NL$^*$ and RNQ are with the Korteweg-de Vries Institute for Mathematics, University of Amsterdam, Amsterdam, The Netherlands. RNQ is also with Centrum Wiskunde \& Informatica, Amsterdam, The Netherlands.} \\
$^*$ corresponding author ({\tt n.a.c.levering@uva.nl})

\vb

\noindent
{\sc Funding:} This research project is partly funded by the NWO Gravitation project N{\sc etworks}, grant number 024.002.003. Date: \today. \\

\end{abstract}

\newpage

\section{Introduction}

The flow of traffic that road networks can carry is not only determined by physical characteristics of roads (width, curvature, inclination, etc.) and traffic rules (e.g., speed limits, prioritization, overtaking), but to a large extent by traffic itself. In the mathematical analysis of road traffic flow this has been recognized in what has been coined the {\em fundamental diagram}:
Large traffic densities may cause traffic speeds to drop rapidly, leading to sudden capacity reductions on individual roads in the network that propagate in space and time, and deteriorate the performance of the entire network.
Most road traffic control mechanisms focus on flow
management inside the network. Such approaches
have proven to be successful in many settings,
but are challenged by the worldwide increase in travel demand.
We therefore consider the real-time \textit{control} of the input traffic streams, 
aiming to prevent the above sketched scenario
in which the factual capacities on individual network roads are exceeded.
By applying such a control procedure, at the cost of a some
waiting time at the on-ramp boundaries of the network, we guarantee
a very small probability of high delays within the network.
Reducing traffic collapses in the network, also limits additional negative consequences
of heavy congestion, such as environmental pollution and economic costs.
The ability to apply control at the boundaries of a road network is
facilitated by advances in Intelligent Transportation Systems. 
Specifically, control may be applied through ramp metering
systems or departure time advice in navigation systems.

When considering traffic flow management, the inherently random nature of
vehicle traffic should be taken into account.
Specifically, due to e.g.\ the heterogeneity in vehicle sizes and 
individual driving habits, the fraction of the capacity of a road that is
unavailable due to the presence of traffic flow on that road
suffers from randomness. 
As the performance of the network,
in terms of realized travel times, is significantly affected 
by catastrophic capacity violations,
(deterministic) control policies that solely consider the 
capacity the {\em average} traffic flow needs on will typically perform poorly.
Thus, there is a need for fast control policies that
do take the fluctuations in traffic flow, and specifically,
random spikes in capacity needs, into account.

In telecommunication networks, similar considerations
for the construction of control policies apply: Transportation mediums have certain {\em bandwidths} (capacities in terms of traffic volume per time),
and each connection using a medium requests part of this bandwidth.
A control policy may decide whether new connections are allowed,
thereby taking into account that the bandwidth requirement of an
individual connection may fluctuate during the time it uses a medium.
Typical performance control focuses on avoiding that the total demand of simultaneous connections exceed the available bandwidth.
A very influential framework in the control of telecommunication
networks is that of \textit{effective bandwidth}, which describes the
minimum bandwidth that needs to be reserved for individual connections
to guarantee a certain level of service for these connections.
The approach leads to linear acceptance
regions, making the effective bandwidth framework a fast and
easily implemented admission control procedure. 
In the context of telecommunications, this framework
was extended to different arrival and service models,
we e.g.\ mention the famous papers by 
\cite{Kelly1991EffectiveQueues,Gibbens1991EffectiveChannel,Elwalid1993EffectiveNetworks}. 
For networks with fluctuating bandwidth demands,
the notion of effective bandwidths was discussed by \cite{Hui1988ResourceNetworks}.
Our goal in this paper is to explore whether the notion of effective bandwidths 
can be extended beyond the telecommunications context,
specifically, to the context of road traffic networks.

{\sc Literature review} \\
There is a broad range of work on the external control of traffic streams in networks, 
in which the control procedures are performed by a traffic planner.
Traditional traffic models, such as the Lighthill-Whitham-Richards 
Model \citep{Lighthill1955OnRoads,Richards1956ShockHighway}, its discretized version, called the Cell Transmission Model~(CTM) \citep{Daganzo1994TheTheory,Daganzo1995TheTraffic}, 
and the Vickrey bottleneck model \citep{Vickrey1969CongestionInvestment,Arnott1990EconomicsBottleneck}, 
operate in a setting in which both demand and delays are deterministic. 
However, in the context of road traffic, uncertainty plays a major role.
Therefore, there is a lot of interest in the stochastic counterparts of these
deterministic flow models, which do account for different driver perceptions, 
moods, car types, etc. 
Recent contributions include the stochastic traffic flow models
of \cite{Jabari2012AFoundations} and \cite{Mandjes2021ANetwork},
the stochastic fundamental diagram of \cite{Qu2017OnApplications},
and the stochastic bottleneck model of \cite{Ghazanfari2021CommuterUncertainty}.

For the routing of individual vehicles, taking uncertainty into
account entails to solving a stochastic shortest path problem,
in which travel times on arcs are (time-dependent) random variables.
Algorithms that minimize the expected travel time or maximize
the on-time arrival probability under various various conditions
are e.g.\ presented in \cite{Levering2021AConditions}. These are stochastic analogues
to (a speed-up versions of) Dijkstra's algorithm, which
yields the optimal route for a vehicle in a deterministic network.
For the optimal routing of traffic \textit{streams} in deterministic
networks, the seminal work of Wardrop introduces
a user equilibrium (i.e., no driver has the
incentive to switch routes) and a social equilibrium
(i.e., minimizing the total network travel cost). 
Examples of stochastic counterparts of
the Wardrop model, in which the
delay is not simply a deterministic function of the traffic flow,
are found in 
\cite{Angelidakis2013StochasticPlayers,Nikolova2011StochasticRouting,Cominetti2019PriceCosts}.
Typically, these studies consider the stochastic user equilibrium,
the stochastic social optimum, or the best or worst
ratio between these as a function of the risk-adverseness
of the vehicles.

Whereas the above works do consider uncertainty, 
their focus lies on the routing of traffic streams.
However, in case of high demand, even with an optimal routing scheme,
unlimited access to a highway network can lead to capacity violations,
which has led to various studies about input rate control strategies.
\cite{Papageorgiou2002FreewayOverview} and 
\cite{Shaaban2016LiteratureMetering} contain overviews of 
ramp metering strategies,
but the referenced works typically consider control in a
deterministic setting. 
In \cite{KovacsThesis}, uncertainty
in the arrival stream is taken into account, but the framework
is limited to a single one-directed road
that consists of multiple segments, and the objective of
study is a proportionally fair control scheme.
A similar objective is studied in \cite{Kelly2010HeavyMotorway},
who do consider uncertainty
in the arrival streams for
a full network of roads. They analyze
the performance of a Brownian network model 
-- often used for proportionally fair control
in telecommunication models -- as
approximate model for the controlled motorway. 

{\sc Contributions} \\
The contributions of this paper are twofold. In the first place, we present
a stochastic flow model, that describes the part of the road capacity that
is effectively taken by the traffic input streams on that road. This overcomes
the limitations of deterministic models, that do not account for
different perceptions, responses, driving habits, car types, etc. 
Specifically, we model the capacity needs as compound Poisson process.
This model offers great flexibility, as we impose only little assumptions on the jump
distribution of this process.

In the second place, we show how the concept of effective bandwidths
can be used to construct a fast control policy in the road traffic
context. This policy allows waiting at the boundaries of the network,
so as to prevent capacity violations within the network.
We show that the asymptotic regime suggests a similar notion
of effective bandwidths, investigate the properties of these effective
bandwidths, and test their optimality through various numerical examples.
In particular, these experiments show that, by applying an effective bandwidth
policy, capacity violations are rare, and that the as a result,
the total time spent before the end of a trip is significantly
smaller than for three other control algorithms that serve
as benchmarks.

{\sc Organization} \\
The remainder of this paper is now structured as follows.
Section~\ref{sec:preliminaries} starts with a description
of the utilized capacity model, to then introduce
the resulting problem of input flow rate control.
The application of effective bandwidths 
for flow control in the context of road
traffic is described in Section~\ref{sec:effband}.
Numerical examples that compare the performance
of our method to three benchmark algorithms are given in
Section~\ref{sec:numericalexamples}.
Section~\ref{sec:conclusiveremarks} contains concluding
remarks.

\section{Preliminaries} \label{sec:preliminaries}

In a road traffic network, 
exceeding the capacity of a road may lead to
extreme delays for (a part of) the road users.
To keep a handle on congestion, we consider
the \textit{control} of the input traffic streams of the network. 
We are, however, challenged by the fact that,
by the heterogeneity within and between traffic input streams,
the impact of traffic flows on the available
capacity suffers from randomness. 
We start by describing a stochastic road traffic model
in Section~\ref{subsec:model}, which describes
the capacity that is effectively used by
the input traffic flows. 
Then, in this model setting, 
Section~\ref{subsec:problem} introduces the problem of
controlling the input flow rates of the network, with the aim
to avert delays due to exceedingly high traffic loads.

\subsection{Utilized capacity model} \label{subsec:model}~\\
We consider a model of road network 
and corresponding graph representation $G = (N,A)$,
of which the set of nodes $N$ represents the ramps
or junctions in the road network and the set of 
directed arcs $A = \{a_1,\dots,a_J\}$ 
represents the roads
connecting these ramps and junctions. 
Each directed arc $a \in A$ has 
a capacity $C_a$, 
indicating how much flow can
be carried by the arc (i.e.,
number of vehicles that
can enter the arc per time unit). 
A path from node $n_1 \in N$
to node $n_2 \in N$ is a collection
of connected directed arcs
that starts at $n_1$ and ends at $n_2$,
and we denote with $\mathcal{P} = \{P_1,\dots,P_I\}$
the set of all paths in the network.
These paths are described by 
the route-link incidence matrix $B$, 
whose elements indicate which links are
part of which paths, i.e.,
\begin{align*}
    B_{i j} =
    \begin{cases}
    1, & a_j \in P_i \\
    0, & a_j \notin P_i
    \end{cases}
    \quad \quad \quad i = 1,\dots,I \quad j = 1,\dots,J.
\end{align*}


Many of the traditional traffic flow models are deterministic.
That is, with $r_{P_i}$ the (mean) input flow rate for
a path $P_i \in \mathcal{P}$
(i.e.,\ the average number of vehicles
per minute traversing $P_i$), it is commonly assumed that
\begin{align*}
    r_{a_j} &= \sum_{i=1}^I B_{ij}r_{P_i}, \quad \quad j  = 1,\dots,J,
\end{align*}
and that the capacity on arc $a_j$ is
exceeded if $r_{a_j} > C_{a_j}$.
However, it is widely recognized that,
due to the variation in
individual driver behaviour and heterogeneity
in vehicle sizes, the capacity needed by the flow
on path $P_i$ is a stochastic measure.
To capture this uncertainty in capacity needs, 
we introduce the random variable $Y_{a_j}$,
denoting the so-called \textit{utilized capacity} of
arc $a_j$, i.e., the flow produced by traffic
traversing arc $a_j$. 
Then, we say that the capacity is
exceeded on arc $a_j$ if the total utilized capacity
is higher than $C_{a_j}$.

For a description of the randomness of the utilized capacity, 
we model $Y_{a_j}$ as a sum of compound Poisson processes. 
We set 
\begin{align}
    Y_{a_j} \myeq \sum_{i=1}^{I} B_{ij} \sum_{k = 1}^{M_i} D_{ik},
    \label{eq:utilizedflow}
\end{align}
with $M_i \sim \text{Poisson}(r_{P_i})$, $r_{P_i} \in \mathbb{R}_{> 0}$,
and $D_{i1},\dots,D_{iM_i}$ i.i.d.
non-negative random variables with known distribution.
Thus, every path $P_i$ that uses arc $a_j$ 
generates a Poisson number $M_i$ of vehicles on that arc, 
with mean $r_{P_i}$, 
the average number of vehicles per minute traversing $P_i$.
The amount that one such vehicle contributes to the occupied capacity
is modeled by $D_{i1} \myeq D_{i}$.
Note that, with limited assumptions
on its distribution, this random variable offers
great modeling flexibility, and can e.g.\ be used to capture
the heterogeneity in vehicle sizes. 
Indeed, modeling the impact of passenger cars 
and trucks with $D_{i}^{\text{cars}}$ and
$D_{i}^{\text{trucks}}$ respectively, 
the impact of the total traffic mix is well described
for $D_{i}$ a mixture distribution of
$D_{i}^{\text{cars}}$ and $D_{i}^{\text{trucks}}$,
the weights set as estimates of their traffic
mix proportions. 
The reader may note the resemblance with
a \textit{passenger car equivalent} factor,
which describes the impact that a mode
of transport has on a traffic variable (in this case, the flow)
compared to a single passenger car 
\citep{Adnan2014PassengerNumbers, Sharma2021EstimationReview}.
Now, as the different paths may contain different traffic mixes,
the occupied capacity is modeled path-dependent.
Concretely, the distribution of $D_{i}$ may differ from
the distribution of $D_{j}$ for $P_i \neq P_j$.



Describing the random vehicular arrivals with a Poisson process
is a natural and widely-used modeling assumption.
The applicability of the Poisson distribution stems from the fact that, 
without congestion or signalized intersections,
drivers behave relatively independent.
In lightly congested traffic conditions,
which form the focus of this study,
empirical observations have indeed shown that Poisson distributed traffic
volumes are realistic.
Note that the Poisson-arrival assumption is also frequently
seen in e.g.\ packet routing and call-center models 
\citep{Koole2013CallOptimization, Bonald2007SchedulingTraffic}. 

\subsection{Avoiding network overflow} \label{subsec:problem}~\\
We consider a network $G = (N,A)$, with input stream rates
$(r_{P_1},\dots,r_{P_I})$, whose impact on the utilized
flow is given through \eqref{eq:utilizedflow}.
It is assumed that, for these given flow demands, traffic in
the network is light, i.e., $\mathbb{P}(Y_a > C_a)$ is sufficiently
small for all $a \in A$. 
Thus, vehicles can travel relatively freely through the network,
experiencing only little hindrance from other road users.
Now, the goal is to decide on the admissibility of additional traffic flow
in a fast and accurate way, such that, with
increased loads, there is an acceptable
balance between the utilized capacity of an arc and the
probability the arc capacity is exceeded, causing congestion.
To this end, we manage potential
increases in the input rates, so as to limit the
number of arc capacity violations.

To avoid reaching the critical capacity, 
we assume we have control of the input stream rates at the boundaries 
of the network, i.e., at the starting nodes
of the paths. Specifically, for each $P_i \in \mathcal{P}$, 
we are able to decide if, instead of $r_{P_i}$, 
a higher input rate 
would still be 
such that capacity violations are
rare, i.e., that $\mathbb{P}(Y_a > C_a)$ is sufficiently
small for all $a \in A$. 
This yields, for all $P \in \mathcal{P}$, a
rule which prescribes whether additional demand can indeed be handled by the
network. If not, the input stream rate of this path should not be
increased, and additional traffic should queue at the boundary of
the network. Our goal is that by controlling traffic in this manner,
at the cost of some delay on the boundaries of the network, extreme congestion within the network, leading to high delays for part
of the traffic, is prevented. 



Now, to construct rules to handle additional demand, we
propose the use of effective bandwidths. 
Effective bandwidths are traditionally applied in 
the management of communication networks,
in which new connections claim part of the available bandwidth,
and bandwidth violations are very undesirable. 
In these networks, effective bandwidths form a powerful tool for admission control:
they efficiently determine a half-space that serves
as acceptance region, such that new connections
are only accepted if they fall into this region. They
are defined arc-wise, as, if congestion is avoided,
the arcs in communication networks are relatively independent
in terms of throughput. Note that, in
vehicle traffic networks, avoiding capacity violations on arcs limits
the negative interaction between the different arcs,
as there are e.g. no traffic jams that affect multiple
arcs. Therefore, there is a promise in expanding
the use of effective bandwidths to the vehicle
network setting, which we explore in this paper.
The application of effective bandwidths
in the context of road traffic will be explained in detail
in the next section.

Before doing so, it is important to remark that applying access control
in vehicular networks is a dynamic procedure.
That is, in case there is additional traffic demand on $P_i \in \mathcal{P}$ that
is allowed into the network, this yields a new input stream rate $\tilde{r}_{P_i}$.
With this new average flow rate, the access rule needs to be updated.
Concretely, for $\varepsilon > 0$ small, it should now be decided
if an input rate of  $\tilde{r}_{P_i}(1 + \varepsilon)$ can still
be handled by the network. Specifically, for the use of effective
bandwidths, this means that, to account for the dynamic updates in traffic streams, 
the acceptance region should be updated regularly.

\section{Effective Bandwidths in Road Traffic} \label{sec:effband}

A concise overview of the concept of
effective bandwidths in their traditional
telecommunications context is provided in
Section~\ref{subsec:background}.
The framework uses a linear acceptance region, 
within which the probability of exceeding the 
network capacity is small, such that the resulting control
procedure is simple and fast. 
Introducing the background
of effective bandwidths, the subsection paves
the way for Section~\ref{subsec:ebub}, in which 
the notion of effective bandwidths is expanded to the
road traffic setting. 


\subsection{Effective bandwidths in telecommunication} \label{subsec:background}~\\
In a telecommunications network, 
many connections are multiplexed over a shared medium. 
The medium has a total bandwidth, and the individual
connections using the medium request part
of this bandwidth. However, typically, the connections
are bursty, in the sense that the bandwidth requirement
may fluctuate during the holding period of the connection.
When applying admission control to such a system,
i.e., when deciding whether a new connection is accepted,
these fluctuations should be taken into account,
as exceeding the bandwidth may lead to connection losses
or other service level violations. Given that each
connection has a mean and a peak rate,
an extreme policy would be to accept a new
connection if the sum of all peak rates does
not exceed the bandwidth, whilst another
extreme policy would be to accept if the sum of all
mean rates is smaller than the bandwidth.
In the first case, there are no service level violations,
but there may be a lot of wasted bandwidth,
given that a connection does not continuously
require its peak rate. In the second case,
the converse is true: there is little excess in bandwidth use,
but there are scenarios in which
service levels are violated.
The concept of effective bandwidths provides a
strategy between these two extremes.

\begin{figure}
    \centering
    \begin{minipage}{0.3\linewidth}
    \centering\captionsetup[subfigure]{justification=centering}
    \includegraphics[width=\textwidth]{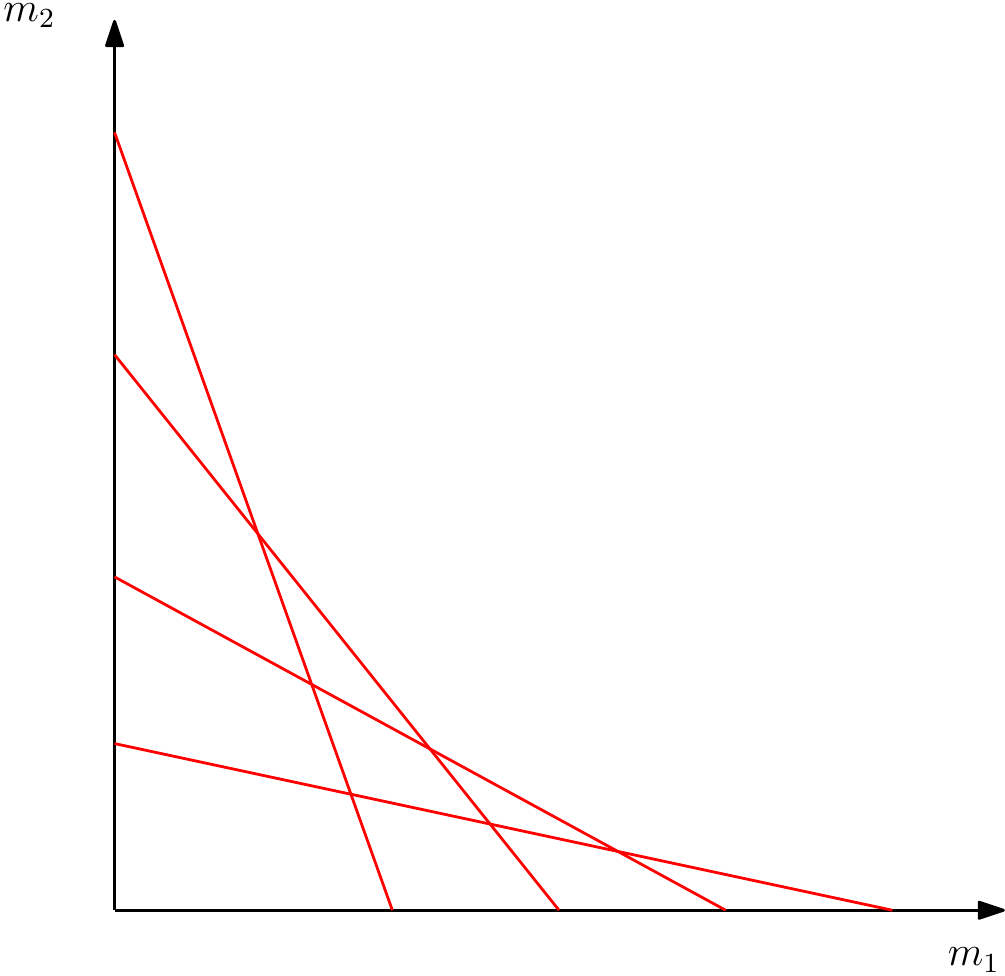}
    \subcaption{}
    \label{fig:detRegionA}
    \end{minipage}
    \begin{minipage}{0.3\linewidth}
    \centering\captionsetup[subfigure]{justification=centering}
    \includegraphics[width=\textwidth]{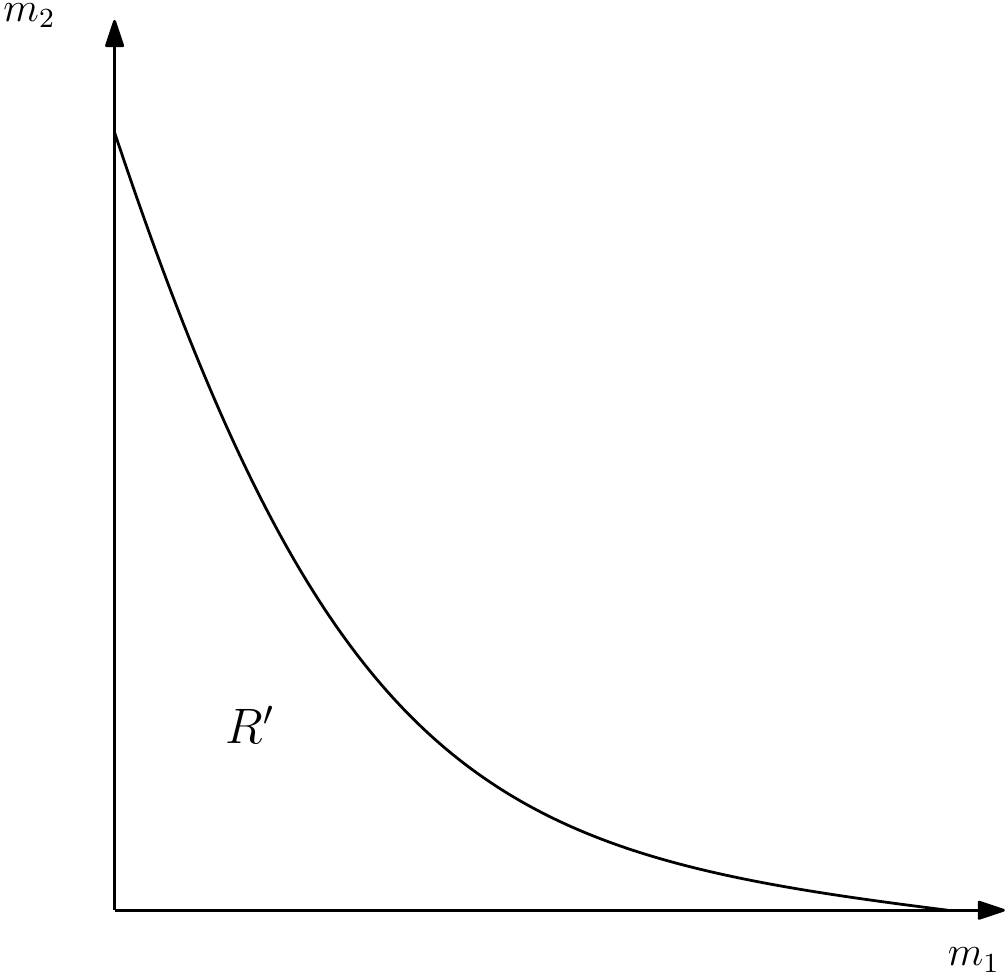}
    \subcaption{}
    \label{fig:detRegionB}
    \end{minipage}
    \begin{minipage}{0.3\linewidth}
    \centering\captionsetup[subfigure]{justification=centering}
    \includegraphics[width=\textwidth]{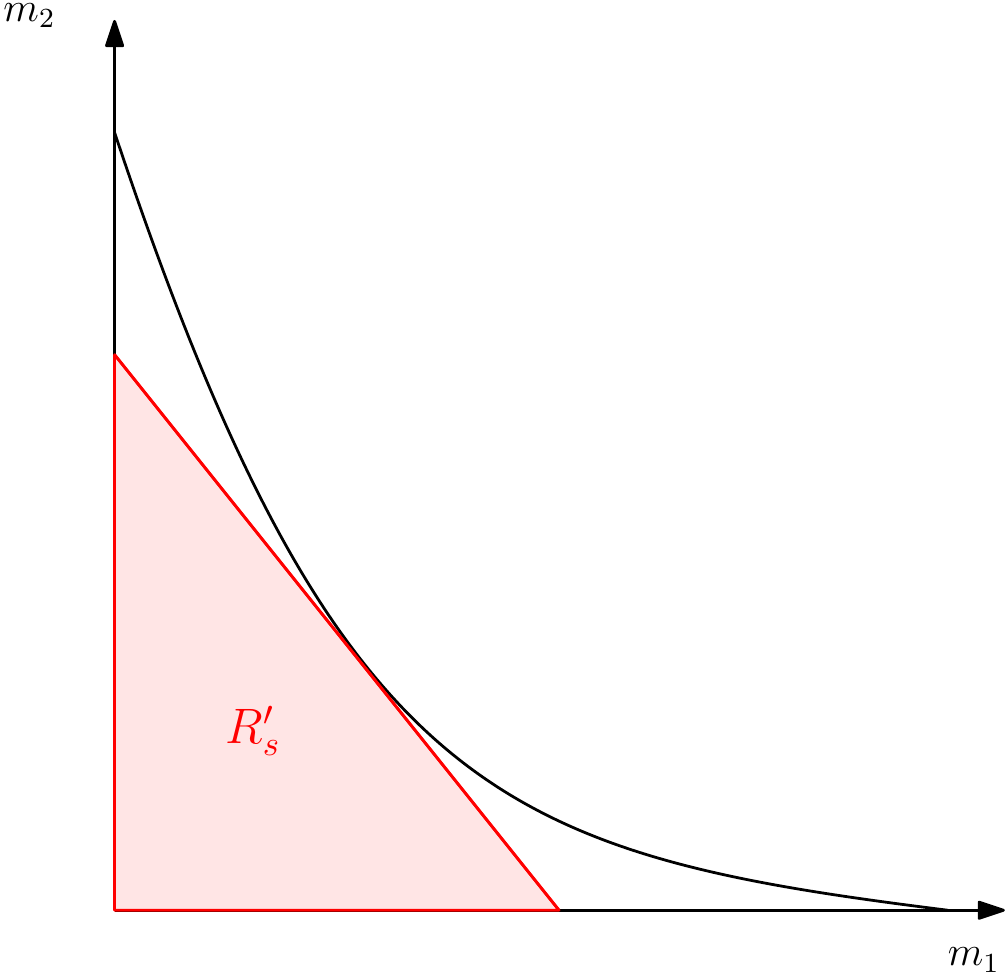}
    \subcaption{}
    \label{fig:detRegionC}
    \end{minipage}
    \caption{Acceptance region and approximation.}
    \label{fig:detRegion}
\end{figure}

Denote $C'$ as the bandwidth shared
by $I'$ types of connections, with $m_i \in \mathbb{N}$ 
connections of type $i = 1,\dots,I'$. Let $D_{ij}'$ be 
the bandwidth requirement for the $j$-th connection of type $i$, 
$D_{i1}',\dots,D_{im_i}'$ identically and independently
distributed, such that
\begin{align*}
    Y' = \sum_{i = 1}^{I'} \sum_{j = 1}^{m_i} \tilde{D}_{ij}'
\end{align*}
is the total bandwidth demanded from the medium.
Observe the similarity between this 
bandwidth demand expression and the road traffic capacity 
demand expression \eqref{eq:utilizedflow}.
Now, parallel to the introduced road traffic setting, the aim
in classical effective bandwidth literature is to 
limit the occurrences in which the
capacity is exceeded. For $\gamma > 0$, the admission control policy
should guarantee $\mathbb{P}(Y' > C') \leq e^{-\gamma}$.
Using a Chernoff bound, it has been proven that
this probability guarantee is satisfied when
the policy is to only accept a new connection
if the new vector of connections falls within the \textit{acceptance region}
$R'$:
\begin{align*}
    R' = \bigcup_{s \geq 0} R_s', 
    \quad 
    R_s' = \Big\{\boldsymbol{m} \in \mathbb{R}_{> 0}^{I'} : 
    \sum_{i = 1}^{I'} m_i \alpha'_i(s) \leq C'-\gamma/s 
    \Big\},
    \quad 
    \alpha'_i(s) = \log \mathbb{E}[\exp{(sD'_{ij})}]/s.
\end{align*}
\begin{remark} \label{rem:cramer}
{\em {
The condition that the new vector should
fall within the region $R'$ is one-way, in 
that satisfying the probability
guarantee does not directly imply that a vector
of connections is within the acceptance region.
However, an application of Cram\'er's theorem shows
that, for large values of $\gamma, C'$
and $m_1,\dots,m_I$, the relation is two-sided.
}} \hfill $\Diamond$
\end{remark}
\begin{example}
{\em {
Observe that the acceptance region $R'$ is a family of half-spaces 
in $\mathbb{R}^{I'}$, indexed by $s$. 
A typical example of the form of $R'$ is
presented in Figure~\ref{fig:detRegion}.
Figure~\ref{fig:detRegionA} shows half-spaces $R'_s$
constructing $R'$ for some values of $s$. The resulting
acceptance region, being the union of these half-spaces, is displayed in Figure~\ref{fig:detRegionB}.
}} \hfill $\Diamond$
\end{example}

Unfortunately, the acceptance region $R'$
is often too difficult to work with. That is,
for a given $I'$-dimensional vector
of connections $\boldsymbol{m}$, it is typically
hard to solve the inversion problem, namely,
to decide whether there \textit{exists} an $s \geq 0$
such that $\boldsymbol{m} \in R'_s$.
Therefore, the idea is to approximate the 
acceptance region $R'$ with 
a region of simpler size.
Specifically, for a chosen $s > 0$, the acceptance
region is approximated with $R'_s$,
whose right boundary is linear (Figure~\ref{fig:detRegionC}).
Then, there is a simple and fast admission control procedure:
only accept a new connection if the new vector
of connections $\boldsymbol{m}$ satisfies
\begin{align*}
    \sum_{i = 1}^{I'} m_i \alpha'_i(s) \leq C'-\gamma/s.
\end{align*}
The weight $\alpha'_i(s)$ of connection type $i$
is called the \textit{effective bandwidth} of type $i$,
and has a value between the mean and peak rate
of the connection type.

\subsection{Effective bandwidths in road traffic} \label{subsec:ebub}~\\
We show how to apply the ideas of the previous subsection
to the road traffic context introduced in Section~\ref{sec:preliminaries}.
Similar to the telecommunications setting, 
the objective is to apply input control to a network, 
so as to bound the probability of exceeding network capacities.
The effective bandwidth framework, working
with a linear acceptance region, provides a fast and 
easily implemented procedure for network admission control.
The effective bandwidths capture both the mean, variance
and other distributional properties 
of the capacity requirements of the different
traffic streams on an arc, and are computed arc-wise,
such that additional traffic on a path is simply
accepted once it satisfies the control constraints
of the arcs making up the path. 


First, with similar techniques as in the
traditional effective bandwidth literature, we
use a Chernoff bound to bound the probability
that the capacity of an arc $a_j \in A$ is exceeded from above:
\begin{equation}
    \log \mathbb{P}(Y_{a_j} > C_{a_j})
    \leq \inf_{s \geq 0} \log \big(\mathbb{E}[e^{sY_{a_j}}]e^{-sC_{a_j}}\big)
    = \inf_{s \geq 0}\big(\sum_{\substack{i = 1,\dots,I \\ B_{ij} = 1}} \log 
    \mathbb{E}[\exp\{s\sum_{k = 1}^{M_i}D_{ik}\}]\!-\!sC_{a_j}\big).
    \label{eq:upperbound}
\end{equation}
Observing that
the expectation in the considered upper bound is the moment-generating function (MGF)
of a compound Poisson distribution, we have, for $\gamma > 0$,
\begin{align}
    \inf_{s \geq 0} \Big(\sum_{\substack{i = 1,\dots,I \\ B_{ij} = 1}} 
    r_{P_i}(\mathbb{E}[e^{sD_{i1}}] - 1)\!-\!sC_{a_j} \Big) \leq -\gamma
    \implies \mathbb{P}(Y_{a_j} > C_{a_j}) \leq e^{-\gamma}.
    \label{eq:ineqEB}
\end{align}
Now, the probability of exceeding the capacity of an arc $a_j$, 
and consequently, the probability
of delay around this arc, is small if $\boldsymbol{r}$ lies in the acceptance
region $R^\alpha_j$:
\begin{align*}
    R^\alpha_j
    = \bigcup_{s \geq 0} R_{j,s}^\alpha, \quad 
    R_{j,s}^\alpha = \Bigg\{\boldsymbol{r} \in \mathbb{R}_{\geq 0}^I : 
    \sum_{i = 1,\dots,I}
    B_{ij} r_{P_i} \alpha_i(s) \leq C_{a_j} - \frac{\gamma}{s}\Bigg\},
    \quad
    \alpha_i(s) = (\mathbb{E}[e^{sD_{i1}}] - 1)/s.
\end{align*}

\begin{remark} 
{\em {
Note that $\alpha_i(s)$ is indeed equivalent to the notion 
of effective bandwidth
as introduced in \cite{Kelly2014StochasticNetworks}.
We observe that,
\begin{align*}
    \alpha_i(s) = \mathbb{E}[e^{sD_{i1}}-1]/s \geq \mathbb{E}[sD_{i1}]/s = \mathbb{E}[D_{i1}],
\end{align*}
such that the bandwidth is between the mean and peak rate of $Y_{a_j}$:
\begin{align*}
    \mathbb{E}[Y_{a_j}] = \sum_{i = 1,\dots,I } B_{ij}r_{P_i} \mathbb{E}[D_{i1}]
    \leq \sum_{i = 1,\dots,I} B_{ij}r_{P_i} \alpha_i(s) < \infty = \sup \{y: \mathbb{P}(Y_{a_j} > y) > 0\}.
\end{align*}
Also, comparing the two different notions of effective bandwidths,
it can easily be seen that $\alpha_i(s) \geq \alpha'_i(s)$. 
This is not surprising, as, 
instead of only protecting against potential high values of $D_{ik}$, 
we need to protect against high values of $M_i$ as well.
}} \hfill $\Diamond$
\end{remark}

\begin{remark} 
{\em {
It is easy to see that \eqref{eq:upperbound} is invariant
to the distribution of $M_i$, in the sense that 
any distribution over the integers would yield a
similar expression for the upper bound.
Notably, 
with $G_X(s) = \mathbb{E}\exp\{sX\}$ the MGF
of a random variable $X$, the implication in \eqref{eq:ineqEB}
may be replaced by the more general
\begin{align*}
    \inf_{s \geq 0} \Big(\sum_{\substack{i = 1,\dots,I \\ B_{ij} = 1}} 
    \log G_{M_i}\big(\log G_{D_{i1}}(s)\big)\!-\!sC_{a_j} 
    \Big) 
    \leq -\gamma
    \implies \mathbb{P}(Y_{a_j} > C_{a_j}) \leq e^{-\gamma},
\end{align*}
such that the results are easily
extended in case $M_i$ is another distribution
from the family of discrete distributions
over the integers for which the MGF is
well-defined.}} \hfill $\Diamond$
\end{remark}

Similar to the telecommunications setting, 
the acceptance region $R^\alpha_j$ is too complex for
practical application, i.e., 
for a given $I$-dimensional vector
of stream rates $\boldsymbol{r}$ it is hard
to decide if there exists an $s_j \geq 0$
such that $\boldsymbol{r} \in R_{j,s_j}^\alpha$.
Therefore, we approximate
the acceptance region 
$R^\alpha_j$ by one of the regions $R_{j,s_j}^\alpha$,
whose right boundary is again linear.
Then, given the regions~$R_{j,s_j}^\alpha$ and
$\varepsilon > 0$, 
the control procedure is simply
to allow an average input rate of $r_{P_k}(1+\varepsilon)$
instead of $r_{P_k}$ if, for all $j \in \{1,\dots,J\}$, the
following inequality is satisfied: 
\begin{equation}
    \sum_{i = 1,\dots,I} B_{ij}r_{P_i} \alpha_i(s_j) 
    + \varepsilon B_{ij} \alpha_k(s_j) \leq C_{a_j}-\gamma/s_j.
    \label{eq:efub}
\end{equation}

For each $j = 1,\dots,J$, we propose
to base the choice of $s_j$ (i.e., $R_{j,s_j}^\alpha$)
on the current traffic conditions.
Given $\boldsymbol{r}$, we let
$R_{j,s_j}^\alpha$ be such that $s_j$
attains the infimum in \eqref{eq:ineqEB}.
With $s_j$ given, the effective bandwidths $\alpha_j(s_j)$
are computed, which may then be stored externally,
such that it can be decided rapidly if \eqref{eq:efub}
is satisfied. Note that, as argued in Section~\ref{subsec:problem},
the application of access control in road networks
is a dynamic procedure. On a longer timescale,
one should therefore adapt the approximating region to new traffic conditions.


\begin{remark}
{\em {
In the above, $\gamma$ is chosen uniformly for all arcs.
We can, however, also work with a bound $\gamma_{a_j} > 0$
per arc $a_j \in A$. This may be preferable if there are 
arcs in the network for which a congestive setting 
has a more deteriorating effect on the network
than congestive settings on other arcs,
as they are e.g.\ located centrally or have a high
degree of neighbouring links.
}} \hfill $\Diamond$
\end{remark}

One drawback of the described 
control method is that \eqref{eq:ineqEB} is not an
equivalence statement,
making the proposed procedure conservative. 
In the classical effective bandwidth notion, the acceptance
region $R'$ is conservative as well, but, in the 
asymptotic regime, $\mathbb{P}(Y' > C') \leq e^{-\gamma}$
is approximately the same as $\boldsymbol{m} \in R'$
(Remark~\ref{rem:cramer}).
Notably, such a limit argument
carries on to the road traffic setting,
as can be deduced from the following theorem.
This theorem can be obtained as a 
consequence of the Bahadur-Rao Theorem
\citep{Bahadur1960OnMean},
but also follows straightforward from the direct argument
presented below.
\begin{theorem} \label{thm:cramer}
For $i = 1,\dots,I$, let $\{M_i(t) : t \geq 0\}$ a Poisson process with
rate $\lambda_i > 0$ and $X_{i1},X_{i2},\dots$ i.i.d.\ random variables
independent of $M_i$ with $\mathbb{E}[X_{i1}] < c_i$,
$\mathbb{P}(X_{i1} > c_i) > 0$, and a log-moment generating function
that is finite for real values in an open neighbourhood of the
origin. Then, with $c = c_1 + \dots + c_I$,
\begin{equation*}
    \lim_{t \to \infty} \frac{1}{t} \log p_t^{c} \equiv
    \lim_{t \to \infty} \frac{1}{t} \log \mathbb{P}\left(\sum_{i = 1}^I X_{i1} 
    + \dots + X_{iM_i(t)}
    > ct \right) = \inf_{s>0}\sum_{i = 1}^I [\lambda_i (\mathbb{E}[e^{sX_{i1}}]-1)-c_is].
\end{equation*}
\end{theorem}
\begin{proof}
The upper bound is an immediate consequence of the Chernoff bound:
\begin{align*}
    \lim_{t \to \infty} \frac{1}{t} \log p_t^c
    &= \lim_{t \to \infty}\frac{1}{t}\log
    \mathbb{P}\left(\sum_{i = 1}^I \sum_{k = 1}^{M_i(t)} X_{ik} > ct\right) 
    \leq \lim_{t \to \infty} \inf_{s > 0} \left\{\frac{1}{t}\log \mathbb{E}\Big[
    \exp \big(s\sum_{i = 1}^I \sum_{k = 1}^{M_i(t)}X_{ik} \big)
    \Big]-cs\right\} \\
    &= \lim_{t \to \infty} \inf_{s > 0} \left\{\frac{1}{t}
     \sum_{i = 1}^I \log
    \mathbb{E}\Big[
    \exp \big(s\sum_{k = 1}^{M_i(t)}X_{ik} \big)
    \Big] -cs \right\}
    = \inf_{s > 0} \sum_{i = 1}^I [\lambda_i \big(
    \mathbb{E}[e^{sX_{i1}}]
    -1 \big) - c_is],
\end{align*}
where the last step follows from \eqref{eq:ineqEB}.
For the lower bound, we note that
\begin{align*}
    \log p_t^c
    &\geq 
    \frac{1}{t} \log \mathbb{P}\left(
    \sum_{i = 1}^I 
    \sum_{k = 1}^{M_i(\lfloor t \rfloor)}
    X_{ik} > c \lfloor t \rfloor \right)
    + \frac{1}{t} \log\mathbb{P}\left(
    \sum_{i = 1}^I 
    \sum_{k = M_i(\lfloor t \rfloor) + 1}^{M_i(t)}
    X_{ik} \geq c(t - \lfloor t \rfloor)
    \right).
\end{align*}
For $M_{i,1}, M_{i,2},\dots$ a sequence
of i.i.d.\ $\text{Poisson}(\lambda_i)$ distributed random
variables,  
\begin{align*}
\lim_{t \to \infty} \frac{1}{t} 
\log \mathbb{P}\left(\sum_{i = 1}^I \sum_{k = 1}^{M_i(\lfloor t \rfloor)}
X_{ik} > c\lfloor t \rfloor \right)
&= \lim_{t \to \infty} \frac{1}{t} 
\log \mathbb{P}\big(\sum_{i = 1}^I \sum_{k = 1}^{
\sum_{l = 1}^{\lfloor t \rfloor} M_{i,l}}
X_{ik} > c\lfloor t \rfloor \big) \\
&\geq \lim_{t \to \infty} \frac{1}{\lfloor t \rfloor} 
\log \mathbb{P}\big(\sum_{l = 1}^{\lfloor t \rfloor} 
\sum_{i = 1}^I \sum_{k = 1}^{M_{i,l}}
X_{ik} > c \lfloor t \rfloor\big) \\
&\equiv \lim_{n \to \infty} \frac{1}{n} \log \mathbb{P}(Z_1 + \dots + Z_n > cn).
\end{align*}
By Cram\'ers theorem, letting $M_i \sim \text{Poisson}(\lambda_i)$ and 
using that $Z_1,\dots,Z_n$ are i.i.d.:
\begin{align*}
    \lim_{n \to \infty} \frac{1}{n} \log \mathbb{P}(Z_1 + \dots + Z_n > cn)
    &= 
    \inf_{s > 0} \left\{ \log \mathbb{E}\big[e^{s \sum_{i = 1}^I \sum_{k = 1}^{M_i} X_{ik}} \big]-cs \right\}
    = \inf_{s > 0} \left\{ \sum_{i = 1}^I [\lambda_i \big(
    \mathbb{E}[e^{sX_{i1}}]-1\big) -c_is] \right\}.
\end{align*}
The theorem now follows from noting that for $t \in \mathbb{R} \setminus \mathbb{N}$
\begin{align*}
    \frac{1}{t} \log\mathbb{P}\left(
    \sum_{i = 1}^I 
    \sum_{k = M_i(\lfloor t \rfloor) + 1}^{M_i(t)}
    X_{ik} > c (t-\lfloor t \rfloor)
    \right)
    &\geq \frac{1}{t} \min_{u \in (0,1]}\log\mathbb{P}\left(
    \sum_{i = 1}^I 
    \sum_{k = 1}^{M_i(u)}
    X_{ik} > cu
    \right) 
    \overset{t \to \infty}{\to} 0.
\end{align*}
\end{proof}

\section{Numerical Experiments} \label{sec:numericalexamples}
Now that we have shown how the notion of effective bandwidths can be used to regulate
traffic streams in road networks, so as to avoid network congestion, we perform
a set of numerical experiments in order to assess the performance of the
proposed policies. Specifically, we demonstrate that capacity violations
are indeed rare, and that
the waiting costs at the boundaries of the network
are typically of a smaller scale than the incurred costs of such violations.
In the experiments, the policy that follows from the effective bandwidth framework is denoted with EB($\gamma$), where $\gamma > 0$ such that the probability of capacity violation is
upper bounded by $e^{-\gamma}$ (see \eqref{eq:ineqEB}).

We compare the performance of the effective bandwidth framework,
in terms of waiting times and number of capacity violations,
to three other input rate control algorithms, which serve as natural
benchmarks. As a first benchmark, we use the procedure that
simply allows all vehicles in the network
at all times, which will be called 'No Control' (abbreviated to NC).
The second benchmark does incur waiting costs at the boundary,
as it will not allow the expected capacity needs on the links
to exceed the link capacities, i.e., $\mathbb{E}[Y_{a_j}] < C_{a_j}$
for all $a_j \in A$.
This procedure will be called
'Expected Needs', and the corresponding policies are abbreviated
with EN. The third benchmark, called 'Random Needs',
takes the random nature of capacity needs into account, and only allows
the current input rate when, for some chosen $\alpha > 0$ and all $a_j \in A$,
\begin{align*}
    \mathbb{E}[Y_{a_j}] + \alpha \sqrt{\text{Var}\left(Y_{a_j}\right)} < C_{a_j}.
\end{align*}
The corresponding procedure is abbreviated to RN($\alpha$).
When aiming for a similar guarantee as in the EB-framework, a natural way
to calibrate $\alpha$ is to use a normal approximation and
set $\alpha \equiv z_{\exp\{-\gamma\}}$, with, for $Z \sim N(0,1)$
and any $\beta > 0$, $z_\beta$ such that $\mathbb{P}(Z \geq z_{\beta}) = \beta$. 
For any $\gamma > 0$, we will denote such a calibration
as $\alpha_\gamma$. 

Experiment~1 examines the impact
of the different policies on a network consisting of a single link,
whereas Experiment~2 considers networks of larger size. 
For the experiments we implemented the networks and
policies in Wolfram Mathematica 12.0 on 
an Intel\textregistered\ Core\texttrademark\ i7-8665U 1.90GHz computer.

\subsection*{Experiment 1}\label{exp:one}
\begin{figure}[h]
    \centering
    \begin{minipage}{0.4\linewidth}
    \centering\captionsetup[subfigure]{justification=centering}
    \includegraphics[width=\textwidth]{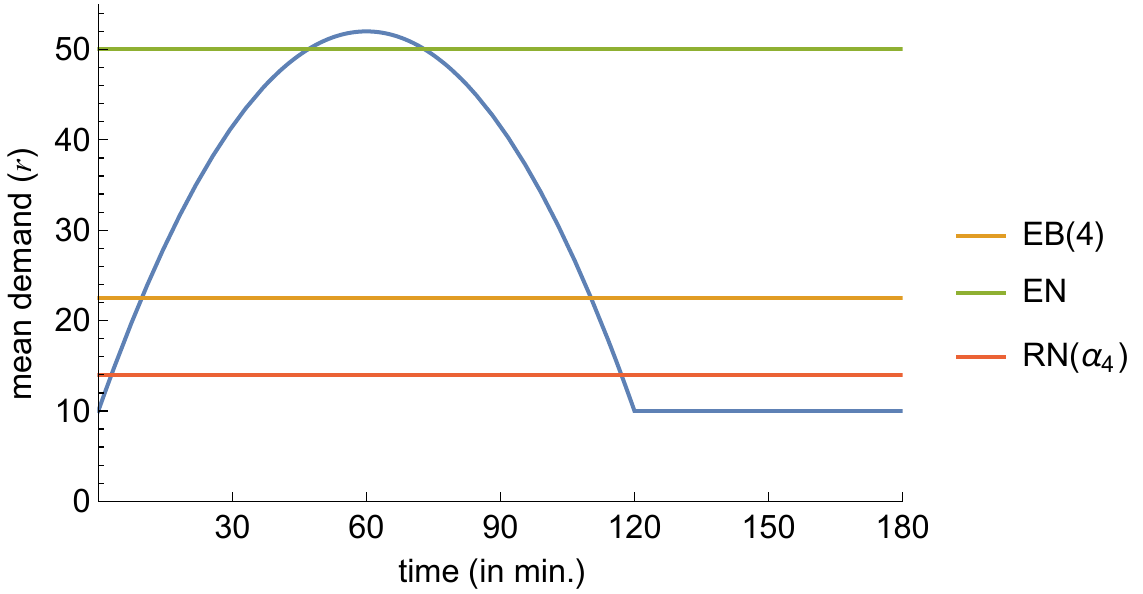}
    \subcaption{Mean traffic demand $r$ as function of time in
    a rush hour setting, together with the maximum input rates of three
    control policies.}
    \label{fig:exp1}
    \end{minipage}
    \qquad 
    \begin{minipage}{0.315\linewidth}
    \centering\captionsetup[subfigure]{justification=centering}
    \includegraphics[width=\textwidth]{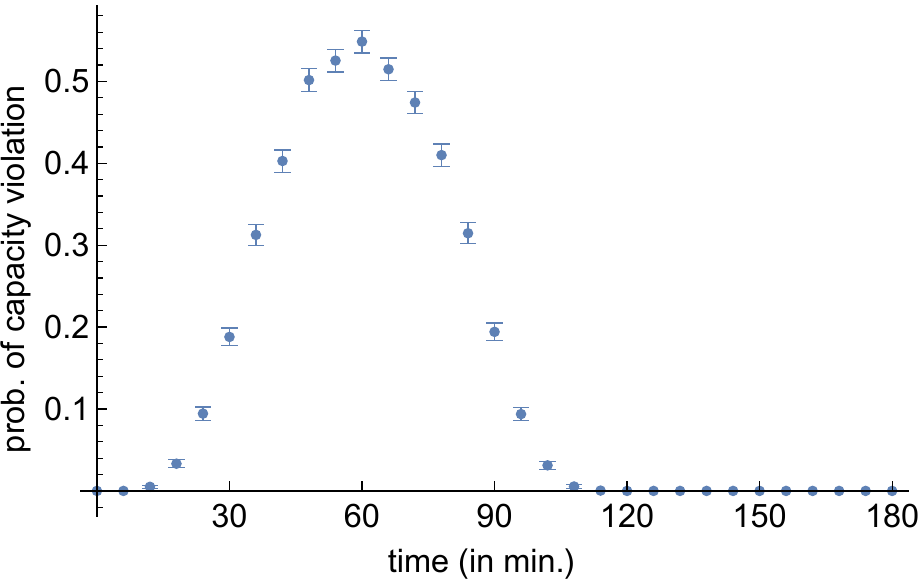}
    \subcaption{Probability of capacity violation as a function of time for policy NC.
    \newline}
    \label{fig:exp1B}
    \end{minipage}
    \caption{Example of the mean demand during rush hour on a single link.}
    \label{fig:exp1Total}
\end{figure}

For illustration purposes, we consider a network that consists of a single link,
with a capacity of 50. 
The amount that one vehicle contributes to the
occupied capacity is modeled as $D \sim$ Hyperexponential($\boldsymbol{p}, \boldsymbol{\lambda}$),
with $\boldsymbol{p} = (0.7,0.3)$ and $\boldsymbol{\lambda} = (3/2,9/16)$. This could
ge.g.\ represent a traffic stream for which any vehicle is a
car  with 70\% probability and with 30\% probability it is a truck; cars and trucks 
occupying a capacity of 2/3 and 16/9, respectively, such that the expected capacity
needs of a single vehicle equal 1.
We consider a rush-hour setting, 
in which, starting at a low mean demand at time 0
(i.e., the time the rush hour initiates), the mean input rate 
of the link first monotonically increases, 
and then monotonically decreases, after which
the rush hour has passed, and the mean
demand stays constant. Specifically, in this experiment,
we consider the mean traffic demand to evolve as in Figure~\ref{fig:exp1}.
For this mean demand curve, Figure~\ref{fig:exp1B} plots 
the (simulated) probability of a capacity violation
at several points in time, 
in case there is no input rate control (policy NC). 
Setting our goal to keep this probability below $e^{-4}$,
we evaluate the network performance for
the policies EN, RN($\alpha_4$) and
EB(4).

Figure~\ref{fig:exp1Total} also 
shows the maximum input rates for the different control policies.
For any of the control policies, the probability of capacity violation can be read off the right hand graph, as long as the mean traffic demand (in the right hand graph) has not yet reached the imposed threshold. For 
the EN-policy this yields a probability
larger than 0.45 that a capacity violation occurs during the peak period of the rush hour.
For the NC-policy, allowing all traffic demand into the network,
this probability is larger than 0.5.
EB(4) and RN($\alpha_4$) are more conservative, and
limit access at an earlier stage. 
Recall that, the latter two policies target at a
probability of capacity violations of $e^{-4} \approx 0.02$ (for RN this is by approximation using a normal distribution, and for EB it is a bound). 


Our aim is 
to avoid the high delays caused by capacity violations, at the expense of some additional delay on the boundaries of the network.
We assess the gains of this approach
in terms of a proxy for the total delay vehicles experience
during their travel, using the following simulation approach.
First, we discretize time in steps
of size $\delta > 0$, and view the evolution
of traffic as a queueing model
that consists of a waiting area of infinite size
at the starting node of the link, 
which we will refer to as the \textit{buffer},
and a FCFS {em queue} on the link itself.
The buffer and queue size at time $t$ are denoted by $B(t)$ and $Q(t)$, respectively.
At each time step $t = \delta m$,
the mean traffic demand $r_t$ is obtained from the area below
the function in Figure~\ref{fig:exp1} between
$\delta(m-1)$ and $\delta m$. The control policy then
decides how much of $B(t\!-\!\delta)$ and $r_t$ is 
offered to the queue.
Let $\bar{\rho}_t$ be the maximum mean flow rate 
the control policy allows at time $t$. 
Then, a mean flow rate of $\bar{r}_t \equiv \min\{\bar{\rho}_t,
B(t\!-\!\delta) + r_t\}$ is admitted to the queue,
such that $B(t) = B(t\!-\!\delta) + r_t - \bar{r}_t$.


Whenever, at time $t = \delta m$, a mean traffic load
of $\bar{r}_t$ enters the link,
the corresponding number of vehicles is simulated,
and they are placed in the queue.
For each of these vehicles, a sample
of their capacity needs yields their service requirement.
The amount of capacity needs that the server
is able to process between $\delta m$ and $\delta (m+1)$ 
is given through a function $c(\cdot)$,
which takes the total capacity needs in the queue 
at time $\delta m$ as input.
The shape of the function is chosen to represent
the impact of high capacity needs on the driveable
speeds. That is, whenever the capacity needs 
are less than the capacity on the link, vehicles
are effectively able to drive the free-flow speed,
represented by the fact that all capacity can be processed. 
Whenever the capacity needs exceed the link capacity,
the attainable speeds are much lower, which
is represented by the fact that $c(\cdot)$ outputs
a value that is far less than the link capacity.
Specifically, in this example, we let
\begin{align} \label{eq:cfunction}
    c(y) \equiv
    \begin{cases}
    y & y \leq 50\delta  \\
    \max\{10\delta,100\delta \!-\!y\} & y > 50\delta . 
    \end{cases}
\end{align}
Note that we multiply the numbers by $\delta$,
to account for the fact that we look at time
windows of size $\delta$.
Furthermore, under congested conditions, we impose
a strictly positive service rate, to make sure
that, for low arrival rates, the model is able to 
recover from these congested conditions
at a future time point.

With the procedure described above we obtain
a proxy for the total delay in the network in the following way.
At each time $t$, an approximation for the total
number of customers in the system is given through the sum
of $B(t)$ and the number of vehicles in the queue.
Computing this sum for each simulation run and each point
in time, we obtain an approximation
for $\tilde{L}$, the average number of customers in the system.
Moreover, averaging over the mean demands $r_t$ of the
different time steps yields an estimate $\tilde{r}$ for the average
arrival rate. Then, with Little's law, we obtain
the proxy $\tilde{L}/\tilde{r}$ for the delay vehicles
in the network experience. The proxies for
different control procedures, 
in the setting described above, with 
the time range of Figure~\ref{fig:exp1},
$B(0) = 0$, $Q(0) = 0$, and $\delta = 1$~min, 
are given in 
Table~\ref{tab:exp1proxies}.
These proxies reveal that the delay for the
NC- and EC-policy is of the same order, which can 
be explained by the fact that EN solely limits access in
the peak of the rush hour. 
The delay under EB(4) is significantly lower, whereas
the delay under RN($\alpha_4$) is of very high order.


\begin{table}[h]
    \centering
    \begin{tabular}{|c|c|c|c|c|c|c|}
    \hline
         & NC & EN & RN($\alpha_4$) & EB(4) \\ \hline
        Delay estimate (in min.) & 49.43 & 49.41 & 66.57 & 41.97
        \\ \hline
    \end{tabular}
    \caption{Estimates for the delays vehicles experience in
    a single-link network under different input rate control
    policies, for the mean demand curve of Figure~\ref{fig:exp1}.}
    \label{tab:exp1proxies}
\end{table}

\begin{figure}[ht]
    \centering
    \includegraphics[width=0.5\textwidth]{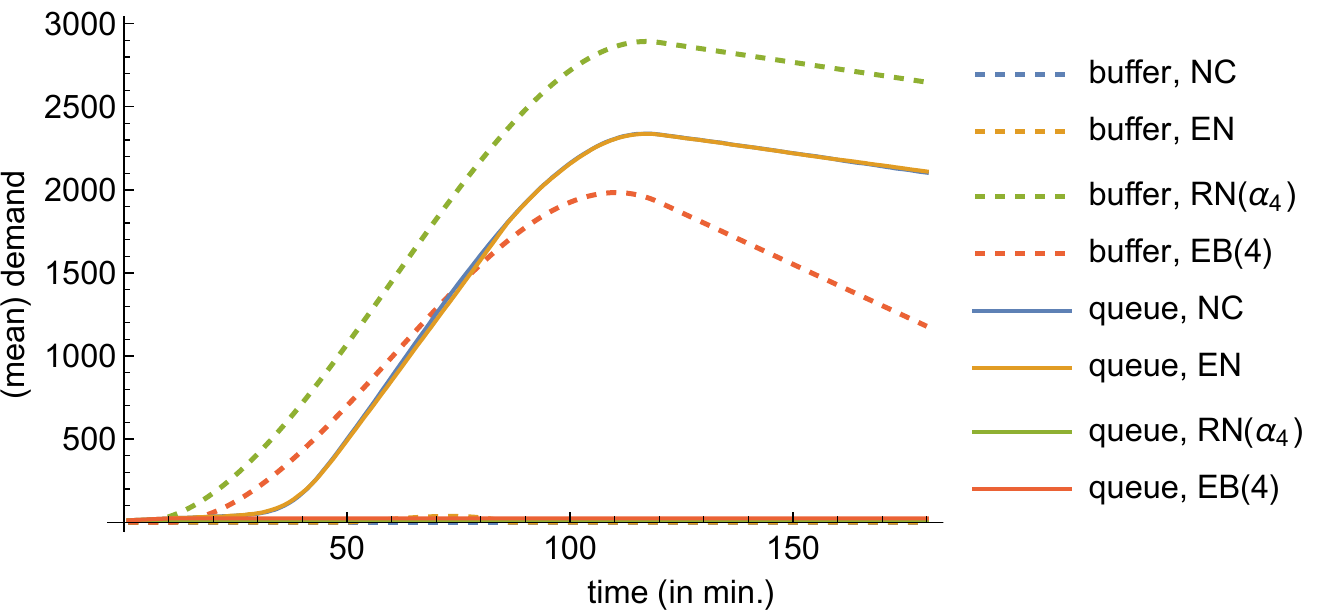}
    \caption{Average buffer and queue size as a function of time, for
    four different control policies.}
    \label{fig:averageSizes}
\end{figure}

The differences in experienced delays as presented in Table~\ref{tab:exp1proxies} 
can also be observed in Figure~\ref{fig:averageSizes}, 
which shows the buffer and average queue size
for the four policies as a function of time.
Note that the mean demand arriving at the buffer is deterministic, such that
the buffer size, in terms of mean demand, 
is a deterministic quantity. The number of vehicles
entering the queue being a random variable, the queue size is a random
quantity, whose average is determined with 10.000 simulations.
Again, we observe that there is little difference between
the NC- and EN-policies, because both allow a high
mean demand onto the link, such that the buffer stays
relatively empty. However, on the link itself, 
due to violations of the link capacity,
the traffic becomes, and stays, highly congested.

Figure~\ref{fig:averageSizes} shows that both EB(4)
and RN($\alpha_4$) indeed succeed in avoiding high
delays on the link, as the queue sizes under both
policies remain of a very small scale. 
Under EB(4), the demand in the buffer grows
during the rush hour period, but as there is no congestion
on the link, the total delay is still much lower
than in the NC- or EN-regime. This is not the case for
the RN($\alpha_4$)-policy: being very conservative, the
buffer size is such that the total delay exceeds
the NC- and EN-regime. However, it is important to
remark that there is a difference between
the delay vehicles experience on the boundaries and
inside the network. That is, by using Intelligent Transportation
Systems, drivers may request their
route from e.g.\ their home, 
and learn when they will be granted
access to the network from there, 
such that waiting at the boundaries
does not automatically yield wasted time for the affected drivers.
This in contrast to waiting inside the network, which,
moreover, results in additional CO$_2$ emissions.





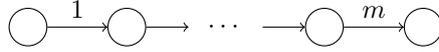
\begin{figure}[ht]
\centering
\begin{tikzpicture}[node distance=1cm,auto]

\node[state, inner sep=5pt,minimum size=5pt, draw=none] (competitive) {};                           
\node (p11a) [state, inner sep=5pt,minimum size=5pt, below of=competitive, yshift=0cm] {};                           
\node (p12a) [state, inner sep=5pt,minimum size=5pt, right of=p11a, xshift=0.3cm] {};           
\node (p13a) [state, inner sep=5pt,minimum size=5pt, right of=p12a, xshift=0.3cm, draw=none] {$\cdots$};
\node (p22a) [state, inner sep=5pt,minimum size=5pt, right of=p13a, xshift=0.3cm] {};
\node (p23a) [state, inner sep=5pt,minimum size=5pt, right of=p22a, xshift=0.3cm] {};

\draw[->] (p11a) -- node[] {$1$} (p12a);
\draw[->] (p12a) -- node[] {} (p13a);
\draw[->] (p13a) -- node[] {} (p22a);
\draw[->] (p22a) -- node[] {$m$} (p23a);

\end{tikzpicture}
\caption{A linear network with $m$ links.} 
\label{fig:linearnetwork}
\end{figure}

\subsection*{Experiment 2}
To examine the impact of capacity violations
in larger networks, we consider the 
linear network of Figure~\ref{fig:linearnetwork},
and vary $m$, the number of links of which the network consists.
The network has a single traffic stream, 
with vehicles wanting
to travel from the left node to the right node.
For this stream, we consider the same mean demand curve
and the same capacity needs characteristics as
in Exp.~1.
Moreover, we let all link capacities equal 100, except for the 
last link, whose capacity is chosen to be 50.

Now, to compute proxies for the delays vehicles
experience under the control procedure,
we discretize time in steps of size $\delta > 0$, and 
again view the evolution
of traffic as a queueing system, with a single buffer
of infinite size at the starting node, and FCFS queues
on each of the $m$ links. The buffer dynamics are 
similar to those in Exp.~1, with the vehicles that are
allowed into the network arriving at the queue that 
corresponds to link~1. These vehicles visit
the queues consecutively: traffic that has been served by the
queue of link~$i$ ($1 \leq i \leq m-1$)
is inserted in the queue of link~$i+1$. If a vehicle
has been served by the queue of link~$m$, it leaves the system.
Now, by summing, for each simulation,
the vehicles in all queues at every time step,
an application of Little's law again yields the order
of travel time delays.

Let $c_i(\cdot)$ be the amount of capacity needs that the 
server at the queue on link $i$ is able to process.
With $\boldsymbol{y}$ 
a vector of the capacity needs on each link, $y_i$
denoting its $i$-th coordinate, and $c(\cdot)$ as in 
\eqref{eq:cfunction}, we 
let $c_m(\boldsymbol{y}) = c(y_m)$.
Then, to capture the fact that traffic jams propagate through
the network, we let, for any link $1 \leq i \leq m-1$, the function $c_i(\cdot)$
explicitly depend on the total capacity needs
in its queue, as well as the total capacity needs in
the next queue. That is, the shape of the function
represents that if the next queue is highly congested,
new traffic is not able to enter this queue, and has to stay in its
current queue. Specifically, allowing not more than $u_i \equiv 190\delta$ per time 
step into queue $i$ for any $1 \leq i \leq m-1$, and not more
than $u_m \equiv 90\delta$ per time step into queue $m$,
\begin{align*}
    c_i(\boldsymbol{y}) \equiv
    \begin{cases}
    \min\{u_{i+1} - y_{i+1} + c_{i+1}(\boldsymbol{y}),y_i\} & y_i \leq 100\delta  \\
    \min\{u_{i+1} - y_{i+1} + c_{i+1}(\boldsymbol{y}),\max\{10\delta,200\delta \!-\!y_i\}\} & y_i > 100\delta . 
    \end{cases}
\end{align*}
With $c_m(\boldsymbol{y})$ known, the amount of vehicles that are
served in each queue can be computed recursively.






\begin{table}[h]
    \centering
    \begin{tabular}{|c|c|c|c|c|c|c|c|c|}
    \hline
         & NC & EN & RN($\alpha_4$) & EB(4) \\ \hline
    $m = 5$ & 49.36 & 48.99 & 68.49 & 45.01 \\
    $m = 10$ & 49.98 & 49.73 & 70.84 & 48.73 \\
    $m = 20$ & 55.43 & 56.52 & 75.30 & 54.81 \\
    $m = 30$ & 62.23 & 61.82 & 79.49 & 60.43 \\
     \hline
    \end{tabular}
    \caption{Estimates for the delays vehicles experience in
    the network of Fig.~\ref{fig:linearnetwork},
    under different input rate control
    policies, for the mean demand curve of Figure~\ref{fig:exp1}.}
    \label{tab:exp2proxies}
\end{table}

The proxies under different input rate control schemes are presented
in Table~\ref{tab:exp2proxies}, and are, naturally, increasing
functions in~$m$. Just as in Exp.~1, the NC- and EN-policy behave relatively
similar. Their high delay values are caused by the fact that the buffer
is relatively empty throughout the complete time window. EB(4) outperforms
these policies, but the high buffer value causes quite high delays, especially
for large values of $m$. However, as argued above, waiting at the boundaries
of the network is significantly different from waiting within the network,
as it is not directly a time waste, and has no negative
environmental consequences. A similar argument can be made for
RN($\alpha_4$), which has the highest delays, but may still be preferred
over NC- and EN.



\section{Concluding remarks} \label{sec:conclusiveremarks}


In this work, we used a compound Poisson process to describe the random part of the
road capacity that is effectively taken by the traffic streams using that road.
To avoid that, in a given road network, link capacities are exceeded, 
we constructed an input rate control policy
that takes the randomness of these capacity needs into account. This policy
guarantees an upper bound on the capacity-violation probability,
and is based on the notion of effective bandwidths as originally introduced in
the telecommunications context. 
Numerical experiments demonstrated that, typically, the total
delay in the network is of a smaller scale than the
delay would be without access control, or with a policy
that only takes expected capacity needs into account.


There are a few natural ways in which our
modeling procedure may be extended, towards which
future work could be specified. For example, in the current
procedure, when
determining the input role control, 
the location of the vehicles
on the paths does not play a role.
Specifically, if given access to the network, 
in our model setting, a demand increase on a given
path will instantaneously lead to an increase
in the capacity needs on all links on this path.
In reality, however, there is a time-component,
and an increase in traffic demand at the boundary will only
lead to an increase in capacity needs on further links
at later points in time.

The aim or our control procedure was to limit the 
negative consequences of capacity violations.
Although our procedure does describe the traffic demand
that is allowed into the network, it does not consider
the fairness in terms of waiting time. That is, if the demand
on a single link consists of two traffic streams 
with the same characteristics, instead of allowing a part
of both streams, our procedure may only allow one of the streams.
Therefore, in terms of practical operationalization,
a potential suggestion for future research would be 
to adapt our procedure so as to meet certain fairness guarantees.


An interesting application of our work would
be the identification of network bottlenecks. 
Since our procedure can identify, 
for a given demand, the set of links whose
capacities have a significant probability
to be violated, this will provide an impression of
the bottlenecks in the network. Notably, 
with the characteristics
of recurrent traffic demand well known, this
may be helpful in deciding on future changes in
network infrastructure.

\section*{Declarations}

\subsection*{Funding}
 This research is funded by the NWO Gravitation project {\sc networks} under grant no.\ 024.002.003.

\subsection*{Competing interests}
The authors have no competing interests to declare that are relevant to the content of this article.

\bibliographystyle{agsm}
\bibliography{references.bib}

\end{document}